\magnification=1200
\line{}
\vskip 1.5 true in
\centerline{\bf Stability of quasi-linear hyperbolic dissipative systems}\par

\bigskip \bigskip
\centerline{by}
\bigskip \bigskip

\centerline{Heinz-Otto Kreiss{\footnote{$^{**}$}{Work supported by
the Office of Naval Research n00014 90 j 1382}} }

\medskip

\centerline{\it  Department of Mathematics, }
\centerline{\it  University of California Los Angeles}
\centerline{\it  Los Angeles CA 90024, USA}

\medskip
\centerline{and}
\medskip

\centerline{Omar E. Ortiz{\footnote{$^{\star}$}{Work partially
supported by CONICET.}} and Oscar A.
Reula{\footnote{$^{\dagger}$}{Researcher of CONICET.}}}

\medskip

\centerline{\it Facultad de Matem\'atica, Astronom\'\i a y F\'\i sica,}
\centerline{\it Universidad Nacional de C\'ordoba,}
\centerline{\it Dr. Medina Allende y Haya de la Torre,}
\centerline{\it (5000) C\'ordoba, Argentina.}


\def\p2{2\pi}

\def\gp2{{1\over \p2}}
\def\partx{{\partial \over\partial x}}

\def\partt{{\partial \over\partial t}}

\def\om{\omega}

\def\calo{{\cal O}}

\def\con{\;{\rm const.}\;}
\def\'{^\prime}

\def\eps{\varepsilon}

\def\snok2#1{\hbox{$\tilde{\mkern -2.0mu \tilde#1}$}}

\bigskip\noindent
{\bf 1. Introduction}

In this work we want to explore the relationship between certain 
eigenvalue condition for the symbols of first order partial differential
operators describing evolution processes and the linear and nonlinear
stability of their stationary solutions.

Consider the initial value problem for the following general first order
quasi-linear system of equations
$$ \displaylines{
v_t=P(v,x,t,\nabla)v =\sum_{\nu=1}^s A_\nu (v,x,t)
{\partial\over \partial x_\nu} v+B(v,x,t)v,\cr
v(x,0)=f(x).\cr}
$$
Here $v$ is a (column) vector valued function of the real space
variables $(x_1,\ldots,x_s)$ and time $t$ with components
$v_1,\ldots,v_n.$ $A_{\nu}$ and $B$ are $n\times n$ matrices and $f(x)$
is a vector valued function of the space variables.

We are interested in solutions which are $\p2$-periodic in all space
variables. There is no difficulty to extend the results to the Cauchy problem
on the whole $x$-space. Instead of Fourier series we would use Fourier
integrals.

We shall restrict our considerations to the case
$$ u_t=\sum_{\nu=0}^s \bigl(A_{0\nu}+\eps A_{1\nu}(x,t,u,\eps)\bigr)
{\partial\over\partial x_\nu} u+
\bigl(B_0+\eps B_1(x,t,u,\eps)\bigr)u.
\eqno(1.1) $$
Here $A_{0\nu},~B_0$ are constant matrices and $\eps$ is a small parameter.
This is, for instance, the case when the stationary solution is constant
and and we consider the solution close to the steady state.
\medskip\noindent
{\bf Assumption 1.1.} {\sl For every $p=0,1,2,\ldots$ and any $c>0,$
there is a constant $K_p$ such that the maximum norm of the $p^{th}$
derivatives of $A_{1\nu},B_1$ with respect to $x,t,\eps$ and $u$ are
bounded by $K_p,$ provided $|u|_\infty\le c.$ For $f(x),$ the
corresponding estimates hold.}
\medskip\noindent
{\bf Definition 1.1.} {\sl The system (1.1) is said to satisfy the
{\bf stability eigenvalue condition}  if there is a constant $\delta >0$
such that, for all real $\om,$ the eigenvalues $\lambda$ of the symbol
$$ \hat P_0(i\om)+B_0:= i\sum_{\nu=1}^s A_{0\nu}\om_\nu +B_0
\eqno(1.2) $$
satisfy
$$ {\rm Re}\, \lambda\le -\delta. \eqno(1.3) $$}
We have to define stability for system (1.1).
\medskip\noindent
{\bf Definition 1.2.} {\sl The system (1.1) is stable if, for any $f,$ 
there exists an $\eps_0$ such that, for $0\le\eps \le \eps_0,$ the
solutions of (1.1) converge to zero for $t\to\infty;$ and there exists an
integer $p_0$ such that $\eps_0$ depends only on the constants $K_p$
with $p\le p_0.$}{\footnote1{We have not specified the norm under
which that convergence takes place, but we shall be using uniform
pointwise convergence.}

In this work we shall look at sufficient conditions under which the stability
eigenvalue condition implies stability.

Consider first the constant coefficient case, i.e., set $\eps =0$ in
the above system. In Section 2 we shall prove that it is possible to
find a positive definite selfadjoint operator $H_0$ such that all
solutions of the system satisfy 
$${d\over dt}(u,H_0u)\le -\delta(u,H_0 u), $$
provided that the problem is well posed in the $L_2$ sense and the
eigenvalue condition is satisfied. In this case the system of equations
is a contraction in a new norm.

In Section 3 we consider linear systems with variable coefficients, i.e.,
the $A_{1\nu}$ depend on $x$ and $t$ but not on $u.$ The construction of
$H$ proceeds via the theory of pseudo-differential operators, i.e.,
we construct the symbol $\hat H(x,t,\om)$ and define the operator $H$ by
$$ 
Hu=\sum_\om e^{i\langle \om,x\rangle}
\hat H(x,t,\om)\hat u(\om)\quad \hbox{for all}\quad 
u=\sum_\om e^{i\langle \om,x\rangle} \hat u(\om).
$$
$\hat H$ depends on the symbols
$$ \hat P_0(i\om)=: \sum_{\nu=1}^s A_{0\nu} i\om_\nu,\quad
\hat P_1(x,t,i\om)=: \sum_{\nu=1}^s A_{1\nu}(x,t) i\om_\nu.
$$
We need that $\hat H$ is a smooth function of all variables. This is only
the case if $\hat P_0,~\hat P_1$ satisfy extra restrictions. 
For the linear and the nonlinear case, we make one of
the following assumptions.
\medskip\noindent
{\bf Assumption 1.2.} {\sl The stability eigenvalue condition is
satisfied and the multiplicities of the eigenvalues of $ \hat
P_0(i\om)+\eps \hat P_1(x,t,u,i\om)$ do not depend on $x,t,u,\om,\eps.$
Also, for every $x,t,u,\om,\eps,$ there is a complete system of
eigenvectors.}
\medskip\noindent
{\bf Assumption 1.3.} {\sl The stability eigenvalue condition is satisfied
and the matrices $A_{0\nu},B_0$ and $A_{1\nu},~\nu=1,\ldots,s,$ are
Hermitian.}

Under any of these conditions we can again construct an $H$-norm and
prove that the problem becomes a contraction.

In the last section we consider the nonlinear equations and the main
result of this paper is
\medskip\noindent
{\bf Main theorem.} {\sl Suppose that Assumption 1.1 and  Assumption
1.2 or 1.3 hold. Then, for sufficiently small $\eps,$ the problem is a
contraction in a suitable $H$-norm and the system (1.1) is thus
stable.}

In the Appendix we relax the eigenvalue condition somewhat.

To prove stability for time dependent partial differential equations
via changing the norm has been done before. For example, in [1] the
method was applied to mixed symmetric hyperbolic-parabolic equations
which included the Navier-Stokes equations. In that case $H$ was
explicitly constructed and not related to an eigenvalue condition. If
we make Assumption 1.3, then our $H$ is similar to the $H$ in [1].

\bigskip\noindent
{\bf 2. Systems with constant coefficients}

In this section we consider the system
$$ \eqalign{
&y_t=\sum_{\nu=0}^s A_{0\nu} {\partial y\over \partial x_\nu} +B_0y=:
\Bigl(P_0(\partx)+B_0\Bigr) y,\cr
&\qquad y(x,0)=f(x),\cr}
\eqno(2.1) $$
with constant coefficients. We are interested in solutions which are
$\p2$-periodic in all space variables. We assume that the problem is well posed
in the $L_2$ sense, i.e., for every $T$ there exists a constant $K(T)$
such that the solutions of (2.1) satisfy the estimate
$$
\|y(\cdot,t)\|\le K(T)\|y(\cdot,0)\|,\quad 0\le t\le T.
\eqno(2.2)
$$
Here
$$ (u,v)=\int_0^{\p2} \cdots \int_0^{\p2} \langle u,v \rangle
dx_1 \cdots dx_s,\quad \|u\|^2=(u,u), $$
denote the usual $L_2$ scalar product and norm. 

One can characterize well posed problems algebraically. Using the Kreiss
matrix theorem (see [2], Sec.2.3), one can prove
\medskip\noindent
{\bf Theorem 2.1.} {\sl The problem (2.1) is well posed in the $L_2$ sense
if and only if it is strongly hyperbolic, i.e., the eigenvalues of the
symbol
$$ \hat P_0(i\om)=i\sum _{\nu=1}^s A_{0\nu} \om_\nu,\quad
\om_j ~{\rm real,} $$
are purely imaginary and, for every fixed $\om\'=\om/|\om|,$ there exists
a complete set of eigenvectors $t_1,\ldots,t_n$ which is uniformly 
independent, i.e., there is a constant $K$ such that
$$ |T^{-1}|+|T|\le K,\quad T=(t_1,\ldots,t_n). $$}
\medskip
We can expand the solution of (2.1) into a Fourier series
$$ y(x,t)=\sum_\om e^{i\langle \om,x\rangle} \hat y(\om,t).
\eqno(2.3) $$
The Fourier coefficients are the solution of the Fourier transformed system 
(2.1)
$$ \hat y_t=\Bigl(i\sum_{\nu=0}^s A_{0\nu}\om_\nu +B_0\Bigr)\hat y=:
\Bigl(\hat P_0(i\om)+B_0\Bigr) \hat y. \eqno(2.4)
$$
We assume that the eigenvalue condition (1.2),(1.3) is satisfied.
Then we can find, for every fixed $\om,$ a positive definite Hermitian matrix
$\hat H,$ a Lyapunov function, such that
$$ \eqalign{
2{\rm Re}\,\hat H(\om)\bigl(\hat P_0(i\om)+B_0\bigr) &=:
\hat H(\om)\bigl(\hat P_0(i\om)+B_0\bigr)
+\bigl(\hat P_0^*(i\om)+B_0^*\bigr)\hat H(\om)\cr
&\le -\delta\hat H(\om).\cr}
\eqno(2.5) $$
Therefore,
$$ \eqalign{
\partt \langle \hat y(\om,t),\hat H(\om)\hat y(\om,t)\rangle 
&=2{\rm Re}\langle \hat y(\om,t),\hat H(\om)
\bigl(\hat P_0(i\om)+B_0\bigr)\hat y(\om,t)\rangle \cr
&\le -\delta 
\langle \hat y(\om,t),\hat H(\om)\hat y(\om,t)\rangle .\cr}
$$
Thus, for every fixed $\om,$ the transformed system (2.4) is a contraction
in the $\hat H(\om)$-norm.

Using the Kreiss matrix theorem, one can prove  (see [2, Sec.2.3])
\medskip\noindent
{\bf Theorem 2.2.} {\sl Assume that the problem (2.1) is well posed in the
$L_2$ sense and that the eigenvalue condition (1.2),(1.3) is satisfied.
Then
we construct the matrices $\hat H(\om)$ such that they
satisfy the uniform inequalities
$$ K_4^{-1}I \le \hat H (\om) \le K_4 I. \eqno(2.6) $$}
Here $K_4$ does not depend on $\om.$
\medskip
We can use $\hat H(\om)$ to define an operator $H$ by
$$ Hu=\sum_\om \hat H(\om)\hat u(\om) e^{i\langle \om,x\rangle}.
\eqno(2.7) $$
It has the following properties

\item{(1)} $H$ is selfadjoint and
$$ K_4^{-1} \|u\|^2 \le (u,Hu)\le K_4 \|u\|^2 . $$
\item{(2)} $ 2{\rm Re}\, H(P_0+B_0)=:
H(P_0+B_0)+(P_0^*+B_0^*)H \le -\delta H.$

These properties follow from Parseval's relation
$$ \eqalign{
(v,Hu) &=\sum_\om \langle \hat v(\om),\hat H(\om)\hat u(\om)\rangle \cr
&=
\sum_\om \langle \hat H(\om)\hat v(\om),\hat u(\om)\rangle =(Hv,u).\cr}
$$
Also,
$$ \eqalign{
K_4^{-1} \|u\|^2 &=K_4^{-1}\sum_\om  |\hat u(\om)|^2
\le \sum_\om \langle \hat u(\om),\hat H(\om)\hat u(\om)\rangle \cr
&=(u,Hu)\le K_4 \|u\|^2\cr}
$$
and
$$ \eqalign{
2\bigl(u,{\rm Re}\, H(P_0+B_0)u\bigr) &=
2\sum_\om \langle \hat u(\om),{\rm Re}\, \hat H(\om)(\hat P_0(i\om)+
B_0)\hat u(\om)\rangle \cr
&\le -\delta \sum_\om \langle \hat u(\om),\hat H(\om)\hat u(\om)\rangle =
-\delta(u,Hu).\cr}
$$
Thus, we can use $H$ to define a new scalar product by
$$ (v,u)_H=(v,Hu),\quad \|u\|^2_H=(u,u)_H,$$
which is equivalent with the $L_2$-norm. The second property gives us
\medskip\noindent
{\bf Theorem 2.3.} {\sl If the conditions of Theorem 2.2 are satisfied,
then the problem (2.1) is a contraction in the $H$-norm.}
\medskip\noindent
{\it Proof.} 
$$ \partt (y,Hy)=2{\rm Re}\bigl(y,H(P_0+B_0)y\bigr)
\le -\delta (y,Hy). $$
This proves the theorem.

\bigskip\noindent
{\bf 3. Linear systems with variable coefficients} 

In this section we want to generalize Theorem 2.2 to linear systems
$$ \eqalign{
v_t &= \sum_{\nu =0}^s \bigl(A_{0\nu}+\eps A_{1\nu}(x,t)\bigr)
{\partial v\over\partial x_\nu} +(B_0+\eps B_1)v\cr
&=:\Bigl(P_0(\partx)+B_0+\eps \bigl(P_1(x,t,\partx)+B_1\bigr)\Bigr)v\cr}
\eqno(3.1) $$
and show that it is a contraction in a 
suitable $H$-norm.
We shall construct
the $H$-norm with help of a pseudo-differential operator
$$ H(t)=H_0+S+\eps H_1(t)  \eqno(3.2)$$
with the following properties.
\smallskip\noindent
(1) $H_0,S,H_1(t)$ are bounded selfadjoint operators. $H_0$ and $S$
do not depend on $t.$ $dH_1/dt$ exists and is also a bounded operator.
Thus, there is a constant $K$ such that 
$$ \|H_0\|+\|S\|+\|H_1(t)\|+\|{dH_1\over dt}\|\le K.$$
\smallskip\noindent
(2) $H_0+S$ is positive definite with
$K$ such that
$$ \|H_0+S\|+\|(H_0+S)^{-1}\|\le K. $$
\smallskip\noindent
$$ 2{\rm Re}\,H_0P_0=:H_0P_0+P_0^*H_0\equiv 0.\leqno(3)$$
\noindent
$$ 2{\rm Re}(H_0+S)(P_0+B_0)=
2 {\rm Re}(SP_0+H_0B_0)\le -\delta (H_0+S).\leqno(4)$$
\smallskip\noindent
(5) $S$ is a smoothing operator with
$$ \|SP_1\|\le K. $$
$$
\|{\rm Re}(H_0+\eps H_1(t))(P_0+\eps P_1)\|=
\eps\|{\rm Re}(H_0P_1+H_1 P_0+\eps H_1P_1)\|\le \eps K. 
\leqno(6)$$
We can prove
\medskip\noindent
{\bf Theorem 3.1.} {\sl Assume that there is an operator $H$ of the form
(3.2) with the properties (1)--(6). For sufficiently small $\eps$ the scalar
product $(u,Hv)$ defines a norm which is equivalent with the $L_2$-norm
and the system (3.1) is a contraction in the $H$-norm.}
\medskip\noindent
{\it Proof.} That $(u,Hv)$ defines a norm which is equivalent with the
$L_2$-norm follows from properties (1) and (2). Also,
$$ \eqalign{
\partt (u,Hu) &=\eps(u,H_{1t} u)+2{\rm Re}\Bigl(u,(H_0+S+\eps H_1)
\bigl(P_0+B_0+\eps(P_1+B_1)\bigr) u\Bigr)\cr
&= \eps (u,H_{1t}u)+
2{\rm Re}\bigl(u,(H_0+S)(P_0+B_0)u\bigr)\cr &\quad +
2{\rm Re}\bigl(u,(H_0+\eps H_1)(P_0+\eps P_1)u\bigr)\cr
&+2\eps {\rm Re}\Bigl(u,\bigl(H_0 B_1+S(P_1+B_1)+H_1(B_0+
\eps B_1)\bigr)u\Bigr)\cr
&\le -\bigl(\delta +\calo(\eps)\bigr) (u,Hu).\cr}
$$
This proves the theorem.
\medskip

We construct the symbol of the pseudo-differential operator (3.2) 
in the following way. Consider all systems with constant coefficients which
we obtain by freezing the coefficients of (3.1) at every point
$x=x_0,~t=t_0.$ We assume that the initial value problem for
all these systems is well posed in the
$L_2$ sense and, therefore, we can construct the matrices $\hat H(x,t,\om)$
for every fixed $x,t.$ Now we think of $\hat H(x,t,\om)$ as a symbol of a
pseudo-differential operator where $x,t$ are independent variables.
Formally, we define the operator $H$ by
$$ Hu=\sum_\om e^{i\langle \om,x\rangle}
\hat H(x,t,\om)\hat u(\om). $$
This definition makes sense only if $\hat H$ satisfies the usual properties of
symbols for pseudo-differential operators. Also, we need the algebra for
such operators to prove that (3.1) becomes a contraction. We want to prove
\medskip\noindent
{\bf Theorem 3.2.} {\sl Assume that the following conditions hold.
\smallskip\noindent
a) There exists a positive definite Hermitian matrix $\tilde H_0(\om\')$ which is a 
smooth function of $\om\'=\om/|\om|$ such that
$$
2{\rm Re}\, \tilde H_0(\om\')\hat P_0(i\om)=:
\tilde H_0(\om\')\hat P_0(i\om)+\hat P_0^*(i\om)\tilde H_0(\om\')
\equiv 0. \eqno(3.3)
$$
\noindent
b) For sufficiently large $|\om|,$ there is a Hermitian matrix
$\tilde S=\tilde S(\om\',1/|\om|)$ which is a smooth function of $\om\'$
and $1/|\om|$ such that
$$
2{\rm Re}\Bigl(\tilde H_0(\om\')+{1\over|\om|} \tilde S(\om\',
1/|\om|)\Bigr)\bigl(|\om|\hat P_0(i\om\')+B_0\bigr)
\le -\delta\bigl(\tilde H_0(\om\')+{1\over|\om|} \tilde S(\om\',1/|\om|)\bigr).
\eqno(3.4)
$$
\noindent
c) There exists a Hermitian matrix $\tilde H_1(x,t,\om\')$ which is a smooth
function of $x,t,\om\'$ such that
$$
2{\rm Re}\Bigl(\tilde H_0(\om\')+\eps \tilde H_1(x,t,\om\',\eps)\Bigr)
\Bigl(\hat P_0(i\om\')+\eps \hat P_1(x,t,i\om\')\Bigr)=0.
\eqno(3.5)
$$

Then we can construct the pseudo-differential operator (3.2) which has
the properties (1)--(6).
Also, there exists an integer $p_0$ such that the constant $K$ depends only
on the first $p_0$ derivatives of the symbols and of the coefficients of
(3.1).
Thus, the problem (3.1) is a contraction in the
$H$-norm.}
\medskip\noindent
{\it Proof.} We construct the symbols for the 
pseudo-differential operators
$$ \eqalign{
H_0 u &=\sum_\om e^{i\langle \om,x\rangle}
\hat H_0(\om)\hat u (\om),\cr
S u &=\sum_\om e^{i\langle \om,x\rangle}
\hat S(\om)\hat u (\om)\cr
H_1 u
& =\sum_\om e^{i\langle \om,x\rangle}\hat H_1 (x,t,\om)\hat u(\om).\cr}
\eqno(3.6)
$$
$\hat H_0(\om),~\hat S(\om)$ do not depend on $x,t.$

Let $C>0$ be a constant. Consider the symbol (1.2) for $|\om|\le C.$
The inequality (1.3) implies (see Lemma 3.2.9 in [2]) that there is
a positive definite Hermitian matrix $\tilde S^{(1)}(\om)$ which is
a smooth function of $\om$ such that
$$
2 {\rm Re}\tilde S^{(1)}(\om)\Bigl(\hat P_0(i\om)+B_0\Bigr)\le
-\delta\tilde S^{(1)}(\om),\quad |\om|\le C+1.
$$
Let $\varphi(|\om|)\in C^\infty$ be a monotone cut-off function
with
$$ \varphi(|\om|)=\cases{ 1 & for $|\om|\ge C+1$\cr
0 & for $|\om|\le C$\cr}.
$$
We define
$$ \eqalign{
&\hat H_0(\om)=\varphi(|\om|)\tilde H_0(\om\'),\cr
&\hat S(\om)={\varphi(|\om|)\over |\om|}
\tilde S(\om\',1/|\om|)+ (1-\varphi(|\om|))\tilde S^{(1)}(\om).
\cr}
$$
It follows from (3.3) and (3.4) that,
for sufficiently large $C,$ the operators $H_0$ and $S$ have
the properties (1)--(5). The symbol
$$ \varphi(|\om|)\Bigl(\tilde H_0(\om\')+
\eps\tilde H_1(\om\',x,t)\Bigr) $$
defines a pseudo-differential operator $H_0+\eps H_{11}$
and the algebra of such operators shows that
$$ H_0+\eps H_{1}=H_0+{\eps\over 2}(H_{11}+H_{11}^*) \eqno(3.7)$$
has the desired properties (1) and (6)
and $K$ can be estimated as required. This proves the theorem.

We shall now give algebraic conditions such that the conditions of
Theorem 3.2 are satisfied.
\medskip\noindent
{\bf Theorem 3.3.} {\sl Assume that Assumption 1.2 holds.
Then we can construct the symbols of Theorem 3.2 whose
derivatives can be estimated in terms of the derivatives of the
coefficients of (3.1).  Therefore, for sufficiently small $\eps,$ the
system (3.1) is a contraction.}
\medskip\noindent
{\it Proof.} We consider the symbol $P_0(i\om) +B_0=|\om|P_0(i\om\')
+B_0$ in a neighborhood of a point $\om\'_0.$ Let $\lambda_1,\ldots,
\lambda_r$ denote the distinct eigenvalues of $P_0(i\om\').$ It is well
known (see, for example [3]) that, because of the constancy of the
multiplicity of the eigenvalues of $P_0(i\om\')$, there exists a smooth
nonsingular transformation $\tilde T_0(\om\')$ such that
$$
\tilde T_0^{-1}(\om\')P_0(i\om\')\tilde T_0(\om\')=
\pmatrix{\Lambda_1 &&0\cr &\ddots &\cr
0 && \Lambda_r\cr}. \eqno(3.8a) 
$$
All eigenvalues of $\Lambda_j$ are equal to $\lambda_j$ and, since
there is a complete set of eigenvectors,
$$ \Lambda_j=\lambda_j I $$
is diagonal.

$\tilde T_0$ is not unique. We can replace it by
$$ T_0=\tilde T_0\pmatrix{ T_{01}&& 0\cr
&\ddots &\cr 0 && T_{0r}\cr}. \eqno(3.8b)
$$
Here the $T_{0j}$ denote arbitrary nonsingular submatrices. We shall
choose them as constant matrices later. (3.8a) gives
$$ \eqalign{
&\tilde T_0^{-1}(|\om|P_0(i\om\') +B_0)\tilde T_0 =: \cr
&\cr
& |\om| \pmatrix{\Lambda_1 &&0\cr &\ddots &\cr
0 && \Lambda_r\cr}+
\pmatrix{\tilde B_{11} & \cdots &\tilde B_{1r}\cr
\vdots&\ddots &\vdots\cr
\tilde B_{r1} &\cdots & \tilde B_{rr}\cr} \cr}
$$
and (3.8b) gives
$$ \eqalign{
& T_0^{-1}(|\om|P_0(i\om\') +B_0) T_0  \cr
&\cr
&= |\om| \pmatrix{\Lambda_1 &&0\cr &\ddots &\cr
0 && \Lambda_r\cr}+
\pmatrix{ T_{01}^{-1}\tilde B_{11} T_{01} & \snok2 B_{12} & \cdots
&\snok2 B_{1r}\cr
&&&\cr
\multispan4 \dotfill \cr
&&&\cr
\snok2 B_{r1} & \snok2 B_{r2} &\cdots &T_{0r}^{-1}\tilde B_{rr}
T_{0r}}.\cr}
\eqno(3.9)
$$
For large $|\om|,$ we can consider the second matrix in (3.9) as
a small perturbation of the first. Therefore, (again, see [2])
there is a smooth transformation $T_1(\om\',1/|\om|)$ such that
$$ \eqalign{
&(I+{1\over|\om|}T_1)^{-1} T_0^{-1}
(|\om|P_0(i\om\') +B_0)T_0
(I+{1\over|\om|}T_1)\cr
&\cr
&= 
|\om| \pmatrix{\Lambda_1 &&0\cr &\ddots &\cr
0 && \Lambda_r\cr}+
\pmatrix{ T_{01}^{-1}\tilde B_{11} T_{01} +
{1\over|\om|}\snok2 B_{11} &&0\cr
&&\cr
&\ddots &\cr
&&\cr
0 && T_{0r}^{-1}\tilde B_{rr} T_{0r} +
{1\over|\om|}\snok2 B_{rr} \cr}.
\cr}
$$
By assumption, the eigenvalues of $\tilde B_{jj}$ have negative real parts.
Therefore, we can choose $T_{0j}$ such that
$$ 2 {\rm Re}(T_{0j}^{-1}\tilde B_{jj} T_{0j})\le -{3 \delta\over 2}
I,\quad |\om\' -\om\'_0| \hbox{ sufficiently small.}
$$
(Again, see Lemma 3.2.9 in [2].) Thus,
$$ \tilde H_0=(T_0^{-1})^* T_0 $$
and, for sufficiently large $|\om|,$
$$
\tilde H_0+{1\over|\om|} \tilde S=\left(\Bigl((I+{1\over|\om|}T_1)T_0\Bigr)^{-1}
\right)^* 
\Bigl((I+{1\over|\om|}T_1)T_0\Bigr)^{-1}
$$
satisfies (3.3) and (3.4). By the usual partition of unity argument, we
can construct $\tilde H_0$ and $\tilde S$ for all $\om\'$ and
conditions (a) and (b) in Theorem 3.2 hold.

We now consider the matrix (symbol)
$$ \hat P_0(i\om\')+\eps \hat P_1(x,t,i\om\'). \eqno(3.10) $$
As the eigenvalues of (3.10) are purely imaginary and their
multiplicity does not change, we can find a smooth transformation
$T_2(x,t,\om\',\eps)$ such that
$$ \eqalign{
& (I+\eps T_2)^{-1}T_0^{-1}(\hat P_0+\eps\hat P_1)T_0(I+\eps T_2)\cr
&\cr
&=\pmatrix{ \Lambda_1 &&0\cr &\ddots &\cr
0 && \Lambda_r\cr}+
 \eps \pmatrix{ \tilde\Lambda_1 &&0\cr &\ddots &\cr
0 && \tilde\Lambda_r\cr}.\cr}
$$
Here $\tilde \Lambda_j=\tilde \lambda_j I$ and $T_2$ is a smooth function
of all variables.
The matrix
$$
\tilde H_0+\eps\tilde H_1=(T^{-1})^* T^{-1},\quad
T=T_0(I+\eps T_2), 
$$
has the property (3.5) and condition (c) in Theorem 3.2 hold.
Therefore, Theorem 3.3 follows from Theorem 3.2.

We consider now the symmetric systems (3.1), i.e., those satisfying
Assumption 1.3. In this case 
the stability eigenvalue condition, for $\om=0,$ implies
that
$$ {\rm Re} B_0 \le -\delta I, \eqno(3.11) $$
and therefore
$${\rm Re} \bigl(u, (P_0+B_0) u\bigr) \le -\delta (u,u).$$
Thus, we can show that (3.1) is a contraction in the usual $L_2$-norm
($H=I$). In the Appendix we shall relax the eigenvalue condition to
some cases where (3.11) does not hold. Therefore we give here a proof
which does not depend on (3.11).

\medskip\noindent
{\bf Theorem 3.4.} {\sl Assume that the coefficients $A_{0j},~A_{1j},~
j=1,2,\ldots,s,$ and $B_0$ but not necessarily $B_1$ are Hermitian
matrices. Assume also that the eigenvalue condition (1.3) holds. Then,
the results of Theorem 3.3 are valid.}

Before we give a proof of the last theorem, we will prove
\medskip\noindent
{\bf Theorem 3.5.} {\sl Assume that, for sufficiently large $|\om|,$
there is a Hermitian matrix $\tilde H(\om)=I+{1\over|\om|}\tilde S$
where $\tilde S=\tilde S(\om\',1/|\om|)$ is a smooth function of
$\om\'$ and $1/|\om|$ such that
$$
2 {\rm Re} \tilde H(\om) \bigl(|\om|\hat P_0(i\om\')+B_0\bigr)
\le -\delta\tilde H(\om). \eqno(3.11)
$$
Then, for sufficiently small $\eps,$ the system (3.1) is a contraction.}
\medskip\noindent
{\it Proof.} The proof proceeds similarly to the proof of Theorem 3.2. It is
much simpler, because in this case we construct a time independent 
pseudo-differential operator of the form
$$ H=I+S $$
which has the properties of Theorem 3.1.
\medskip\noindent
{\it Proof of Theorem 3.4.} Consider the symbol
$|\om|\hat P_0(i\om\')+B_0$ for large $|\om|.$ Let $\om\'=\om_0\'$ be
fixed. Since the coefficients $A_{0j}$ are Hermitian, there is a unitary
transformation such that
 $$ \eqalign{
 & U^*(\om'_0)\Bigl(|\om|\hat P_0(i\om'_0)+B_0\Bigr) U(\om'_0)\cr
 &= i|\om|\pmatrix{ 
 \Lambda_1 & & 0\cr &\ddots &\cr 0&& \Lambda_r}+
 \pmatrix{\tilde B_{11} &\tilde B_{12} & \cdots & \tilde B_{1r}\cr
 \tilde B^*_{12} &\tilde B_{22} &\cdots &\tilde B_{2r}\cr
 \vdots &&\ddots &\vdots \cr
 \tilde B^*_{1r} &\cdots &\tilde B^*_{r-1\, r} &\tilde B_{rr}\cr} \cr}
 \eqno(3.12)
 $$
Here $$ \Lambda_j=\lambda_j I$$
represent the different eigenvalues according to their multiplicity.
Since $\tilde B_{jj}$ are also Hermitian, we can assume that they are
diagonal.
Otherwise, we apply a block-diagonal unitary transformation to (3.12).
For large $|\om|,$ we consider the $B$-matrix in (3.12) a small
perturbation of $i|\om|\Lambda.$
Therefore, we can construct a transformation $I+{1\over |\om|} T(\om'_0)$
such that
$$ \eqalign{
& \Bigl(I+{1\over|\om|}T(\om'_0)\Bigr)^{-1} U^*(\om'_0)
 \Bigl(|\om| \hat P_0(i\om'_0)+B_0\Bigr)U(\om'_0)
 \Bigl(I+{1\over|\om|}T(\om'_0)\Bigr)\cr
 &= i|\om|\pmatrix{ 
 \Lambda_1 & & 0\cr &\ddots &\cr 0&& \Lambda_r}+
\pmatrix{\tilde B_{11} & & 0\cr
&\ddots &\cr
0 && \tilde B_{rr}\cr} +{1\over |\om|}\snok2 B\cr
&=:i|\om|\Lambda +\tilde B+{1\over |\om|}\snok2 B.\cr}
$$
The eigenvalue condition guarantees that $\tilde B_{jj}\le -\delta I$
for all $j$ and, for sufficiently large $|\om|,$
$$
2{\rm Re}(i|\om|\Lambda +\tilde B+{1\over|\om|}\snok2 B)\le
-{3\over 2}\delta I. $$
We shall now show that there is a neighborhood of $\om\'_0$ where the
matrix $\tilde H(\om)$ of (3.11) is given by
$$
\tilde H(\om)=U(\om_0\')\bigl(I+{1\over|\om|}T^*(\om_0\')\bigr)^{-1}
\bigl(I+{1\over|\om|}T(\om_0\')\bigr)^{-1}U^*(\om_0\')
=:I+{1\over|\om|} \tilde S(\om_0\',{1\over |\om|}).
$$
We have
$$ \eqalign{
2{\rm Re}\,& \tilde H(\om)\bigl(|\om|P_0(i\om\')+B_0\bigr)\cr
&=
2{\rm Re}\,\tilde H(\om\')\bigl(|\om|P_0(i\om_0\')+B_0\bigr)+
|\om|
2{\rm Re}\,\tilde H(\om)P_0\bigl(i(\om\'-\om_0\')\bigr)\cr
&\le -{3\over 2}\delta \tilde H(\om)+|\om|\cdot
{1\over|\om|}2{\rm Re}\, \tilde S(\om_0\',{1\over|\om|})P_0\bigl(i(\om\'-\om_0\')\bigr)\cr
&\le \Bigl(-{3\over 2} \delta +\con |\om\'-\om_0\'|\Bigr) \tilde H(\om).}
$$
Thus, for sufficiently small $|\om\' -\om_0\'|,$ the inequality (3.11) holds.
With help of the usual partition of unity argument (see again Lemma 3.2.9
of [2]), we can construct $\tilde H(\om)$ for all $\om\'$ and the theorem
follows from Theorem 3.5.

\bigskip\noindent
{\bf 4. Nonlinear systems.} 

In this section we consider the nonlinear system (1.1). We start with
the case that $A_{0\nu},A_{1\nu},~\nu=1,\ldots,s;$ are Hermitian
matrices and
$${\rm Re}\, B_0\le -\delta. \eqno(4.1) $$
Our arguments follow closely the arguments in [2, Chapter 5,6] and we 
assume that the readers are familiar with them.

\noindent
We shall derive a priori estimates and shall use the following notations:
$j=(j_1,\ldots,j_s),~j_\nu$ natural numbers, denotes a multi-index,
$|j|=\sum j_\nu,~D^j=\partial ^{j_1}/\partial x_1^{j_1} \cdots 
\partial ^{j_s}/\partial x_s^{j_s} $ denote the space derivatives and
$$ \|u\|^2_p=\sum_{|j|\le p} \|D^j u\|^2 $$
denotes the derivative norm of order $p.$

To begin with, we assume that $\eps =0$ and derive estimates for
$$ \eqalign{
&{\partial u\over\partial t}=\bigl(P_0(\partx) +B_0\bigr)u,\cr
&u(x,0) = f(x).\cr} \eqno(4.2) $$
Differentiating (4.2) gives us
$$ (D^j u)_t=P_0(\partx) D^j u +B_0 D^j u.$$
Therefore, by (4.1),
$$ \eqalign{
\partt \|D^j u\|^2 &= 2{\rm Re}\bigl(D^j u,P_0(\partx)D^j u\bigr)
+ 2{\rm Re}\bigl(D^j u,B_0 D^j u\bigr)\cr
&= 2{\rm Re}\bigl(D^j u,B_0 D^j u\bigr)\le -2\delta \|D^j u\|^2.\cr}
$$
Adding these inequalities for all $j$ with $|j|\le p$
we obtain, for any $p,$
$$ \partt \|u\|^2_p \le -2\delta \|u\|^2_p, \eqno(4.3) $$
i.e.,
$$ \|u(\cdot,t)\|^2_p \le e^{-2\delta t} \|u(\cdot,0)\|^2_p. $$
Now we consider the nonlinear system (1.1). We derive an estimate for 
$p\ge s+2.$ Local existence causes no difficulty, it has been known for
a long time. There exists an interval $0\le t\le T,~T>0,$ where the
solution exists and
$$ \|u(\cdot,t)\|^2_p \le 2 \|u(\cdot,0)\|^2_p. \eqno(4.4)$$
There are two possibilities:
$$ \hbox{Either}\quad T=\infty\quad \hbox{or}\quad
T<\infty \hbox{ and } \|u(\cdot,T)\|^2_p=
2\|u(\cdot,0)\|^2_p. \eqno(4.5) $$
We shall now prove that $T=\infty$ for sufficiently small $\eps_0$ and
that the initial value problem is a contraction (see [2, Section 6.4.1])

We differentiate (1.1) and obtain
$$
(D^j u)_t=\Bigl(P_0(\partx)+B_0\Bigr) (D^j u)+
\eps P_1(x,t,u,\partx) (D^j u)+\eps R_j,
$$
where $R_j$ denote lower order terms. Therefore,
$$ \eqalign{
\|D^j u\|^2_t &= 2{\rm Re}\Bigl(D^j u,\Bigl(P_0(\partx)+B_0\Bigr) 
D^j u\Bigr)\cr
&\quad +2{\rm Re}\,\eps \Bigl(D^j u, P_1(x,t,u,\partx)D^j u\Bigr)+
2\eps {\rm Re}(D^j u,R_j).\cr} \eqno(4.6)
$$
By (4.1),
$$
2{\rm Re}\Bigl(D^j u,\Bigl(P_0(\partx)+B_0\Bigr)  D^j u\Bigr)
\le -2\delta \|D^j u\|^2. \eqno(4.7)
$$
Integration by parts gives us
$$ \eqalign{
{\rm Re} \Bigl(D^j u, P_1(x,t,u,\partx)D^j u\Bigr)&=
-{1 \over 2} \sum_{\nu=1}^s \Bigl(D^j u,{\partial A_{1\nu}\over\partial x_\nu} D^j u \Bigr)\cr
&\le \con \; K_1 \Bigl(1+ \sum_{\nu=1}^s \left|{\partial u\over
\partial x^{\nu}}\right|_\infty \Bigr) \|D^j u\|^2\cr
& \le M_1 \|u\|^2_p.\cr} \eqno(4.8)
$$
The $M_j$ are polynomials in $\|u\|_p$ of degree $|j|$ whose coefficients
depend only on the constants $K_0,\ldots, K_j$ of Assumption 1.1. Using
Sobolev inequalities we find bounds
$$ \|(D^j u,R_j)\| \le M_j \|u\|_j^2. $$
Adding all these inequalities gives us
$$
\partt \|u\|^2_p \le -2\delta \|u\|^2_p+\eps M \|u\|_p^2,\quad
M=\max_j M_j. \eqno(4.9)
$$
Thus, for all $\eps$ with $0\le \eps\le \eps_0$ with $\eps_0$ 
sufficiently small, we have
$$ \partt \|u\|^2_p \le -\delta \|u\|^2_p. $$
Therefore, $T=\infty$ and the initial value problem is a contraction.
\medskip
We now consider the general case. We assume that the assumptions of Theorem
3.3 or 3.4 are satisfied. Again, we begin with the case that $\eps=0.$ 
Then there is a pseudo-differential operator
$$\tilde H=H_0+S $$
which defines a norm that is equivalent with the $L_2$-norm such that
$$
\partt \|u\|^2_{\tilde H}=\partt (u,\tilde Hu)=
2{\rm Re}(u,\tilde H\Bigl(P_0+B_0)u\Bigr)\le -\delta (u,u)_{\tilde H}.
$$
Thus,
$$ \|u(\cdot,t)\|^2_{\tilde H}\le e^{-\delta t} \|u(\cdot,0)\|^2_{\tilde H} $$
and, therefore, also
$$
\|u(\cdot,t)\|^2_{\tilde H,p} \le
e^{-\delta t}\|u(\cdot,0)\|^2_{\tilde H,p},\quad
\|u\|^2_{\tilde H,p}=\sum_{|j|\le p} \|D^j u\|^2_{\tilde H}.
$$

Now we consider the general system (1.1). We proceed
in the same way as in the previous case and derive estimates for
$p\ge s+2.$ The only difference is that we derive the estimates in
the $H$-norm.

Local existence is again no difficulty. Thus, there is an interval
$0\le t\le T,~T>0$ where
$$ \|u(\cdot,t)\|^2_{\tilde H,p} \le 2\|u(\cdot,0)\|^2_{\tilde H,p}. $$
For $T,$ the alternative (4.5) holds. In this interval we can estimate
the solution and its derivatives up to order $p-[s/2]-1$ in the maximum
norm in terms of $\|u(\cdot,0)\|^2_{\tilde H,p}.$ Thus,
we can think of the system (1.1) as a linear
system and construct the pseudo-differential operator (3.2) 
and estimate the solution and its derivatives in the $H$-norm which
differs from the $\tilde H$-norm only by terms of order $\eps.$ The symbol
depends on the solution but, by Theorem 3.2, if $p$ is sufficiently large,
then the constant $K$ in Theorem 3.1 can also be estimated in terms of
$\|u(\cdot,0)\|^2_{\tilde H,p}.$ The rest of the proof proceeds as before. We
differentiate (1.1) with respect to the space derivatives and obtain
in the $H$-norm
$$ \eqalign{
\partt (D^j u,HD^j u)&=\eps(u,H_{1t} u)+
\Bigl(D^j u,H\bigl((P_0+B_0+\eps(P_1+B_1)\bigr) D^j u\Bigr)\cr
&\quad +2\eps {\rm Re}(D^j u,H R_j).\cr}
$$
Using Theorem 3.1, we obtain the inequality (4.9) but now in the $H$-norm.
Thus, we have proved the Main Theorem of Section 1.

\bigskip\noindent
{\bf Appendix}

We want to relax the stability eigenvalue condition for the cases when
some of the eigenvalues of $B_0$ have zero real part. We do this only
for the 2$\pi$-periodic case.
\medskip\noindent
{\bf Definition A.1.} {\sl The system (1.1) is said to satisfy the {\bf
relaxed stability eigenvalue condition}  if the following conditions
hold.

\noindent
1) There is a constant $\delta >0$ such that, the eigenvalues
$\lambda$ of the symbol $\hat{P}(i\omega) + B_0$ satisfy
$$ {\rm Re}\, \lambda\le -\delta \eqno(A.1) $$
for all $\om = (\om_1, \ldots, \om_s) \neq 0,$ $\om_j$ integer. 

\noindent
2) The eigenvalues $\lambda(0)$ of $B_0$ satisfy
$$ \hbox{Either $\quad$ Re}\,\lambda \le -\delta \quad
\hbox{or}\quad \lambda =0. \eqno(A.2)
$$
Also, if the multiplicity of the zero eigenvalue is $r,$ then there are
$r$ linearly independent eigenvectors connected with $\lambda=0.$

\noindent
3) The nullsapce of $B_1$ contains the nullspace of $B_0.$ }
\medskip
We can find a nonsingular transformation $S$ such that
$$ S^{-1}B_0 S=\pmatrix{B_{01} & 0\cr 0 & 0\cr},
\quad B_{01} \hbox{ nonsingular.} \eqno(A.3)
$$
If $B_0$ is symmetric, we can choose $S$ to be unitary. Therefore, we can
assume that $B_0$ already has the form (A.3). Then, by the third part of the
assumption, $B_1$ has the form
$$ B_1=\pmatrix{ B_{11} & 0\cr B_{12} & 0\cr}. \eqno(A.4) $$
Let
$$ u(x,t)=\pmatrix{ \hat u^{I} (0,t) \cr \hat u^{II} (0,t) \cr}
+\sum_{\om\ne 0} e^{i\langle \om,x\rangle} \hat u(\om,t).
$$
Here the partition of $\hat u(0,t)$ corresponds to that of $B_0,B_1.$
Denote by $Q$ the projection
$$ Q u(x,t)=\pmatrix{0\cr \hat u^{II} (0,t) \cr}. $$
Using the notation
$$ Q u(x,t)=:u^{(0)} (t),\quad
(I-Q) u(x,t)=:v(x,t),$$
we can write the system (1.1) as
$$ \eqalign{
u_t^{(0)} &= Q(P_0+B_0)(u^{(0)} +v)+\eps 
Q(P_1+B_1)(u^{(0)} +v)\cr
&=\eps Q(P_1+B_1)v.\cr} \eqno(A.5)
$$
$$ \eqalign{
v_t &= (I-Q)(P_0+B_0)(u^{(0)} +v)+\eps 
(I-Q)(P_1+B_1)(u^{(0)} +v)\cr
&=(P_0+B_0)v+\eps (I-Q)(P_1+B_1)v.\cr} \eqno(A.6)
$$
As before, we need only to consider linear systems. Then (A.6) decouples
completely from (A.5). It is a system on the subspace $(I-Q)L_2.$
Our results tell us that, for sufficiently small $\eps,$
it is a contraction and $v$ converges exponentially to zero. Since
$$ u^{(0)} (t)=\eps \int_0^t Q(P_1+B_1)v(x,\xi) d\xi + u^{(0)}(0), $$
it follows that also $u^{(0)}(t)$ converges for $t\to\infty.$

We summarize the results of the appendix in the following theorem.
\medskip\noindent
{\bf Theorem A.1.} {\sl Suppose that assumption 1.1 and assumption 1.2
or assumption 1.3 hold but with the stability eigenvalue condition
replaced by the relaxed stability eigenvalue condition. Then, for
sufficiently small $\eps$, the problem is a contraction, in a suitable
$H$-norm, for the nontrivial part $v$ of the solution of (1.1) and
$u^{(0)} \rightarrow \con$ when $t\rightarrow \infty$. Thus, the system
(1.1) is stable in this generalized sense.}


\bigskip\noindent
{\bf References}

\bigskip

\item{[1]} T. Hagstrom and J.Lorenz ,``All-time existence of smooth
solutions to PDEs of mixed type and the invariant subspace of uniform
states", {\it Advances in Appl. Math.}, 16, pp. 219-257, (1995).

\item{[2]} H.O. Kreiss and Jens Lorenz, {\it Initial-Boundary Value
Problems and the Navier-Stokes Equations}, Academic Press, (1989).

\item{[3]} T. Kato, {\it Perturbation Theory for Linear Operators,}
Springer Verlag, (1980). See also, G.W. Stewart and Ji-guang Sun, {\it
Matrix Perturbation Theory,} Computer Science and Scientific
Computing, Academic press Inc., (1990).

\bye